\documentclass[10pt]{article}
\usepackage{amsmath,amssymb}
\usepackage{latexsym, amsfonts, layout, fancyhdr, setspace}

\newtheorem{theo}{Theorem}[section]

\newtheorem{lem}[theo]{Lemma}

\newtheorem{df}[theo]{Definition}
\newtheorem{rk}[theo]{Remark}

\newcommand{\R}{\mathbb R}

\newcommand{\B}{\mathcal{B}}
\newcommand{\U}{\mathcal{U}}
\newcommand{\EL}{\mathcal{L}}

\makeatletter \oddsidemargin-.25in \evensidemargin-.25in
\makeatother

\begin{document}

\date{}
\title{\Large\bf Global attractors for the damped nonlinear wave equation in unbounded
domains }
\author{Djiby Fall\\[2mm]
 Mathematics Program, Division of Science\\ NewYork University Abu Dhabi, UAE\\
Email: dfall@nyu.edu, falldjiby@gmail.com\\[2mm] \&\\
Yuncheng You\\[2mm]
Department of Mathematics and Statistics\\  University of South Florida\\
Email: you@mail.usf.edu\\[2mm]}

\maketitle

{\footnotesize \noindent {\bf Abstract.}  The existence of a global attractor for wave equations in unbounded domains is a challenging problem due to the non-compactness of the Sobolev embeddings. To overcome this difficulty, some authors have worked with weighted Sobolev spaces which restrict the choice of the initial data. 
Using the "tail estimation method" introduced by B. Wang for reaction diffusion equations, we establish in this paper the existence of a global attractor for two wave equations in the traditional Hilbert spaces $\displaystyle H^1(\Omega)\times L^2(\Omega)$ where $\Omega$ is an unbounded domain of $\R^N$. The first equation, with a mass term is studied in the whole space $\R^N$ and the second one without mass term is considered in a domain bounded in only one direction so that Poincar\'e inequality will hold.
\\\\
{\bf Keywords.} wave equation, global attractor, unbounded domain, asymptotic compactness, absorbing set, ``tail estimates'', Poincar\'e Inequality.\\

\vskip.2in

\section{Introduction}
We study in this paper the existence of a global attractor for the
following two damped nonlinear wave equations in an unbounded domain
of $\R^N$:
\begin{eqnarray}\label{inwav1}
u_{tt} + \lambda u_t-\triangle u+u+f(u)=g(x),  \quad t>0
\end{eqnarray}
and
\begin{eqnarray}\label{inwav2}
u_{tt} + \lambda u_t-\triangle u+f(u)=g(x), \quad t>0
\end{eqnarray}
where $\lambda$ is a positive constant, $g$ is a given function
and $f$ is a nonlinear term satisfying some growth conditions to
be specified later. The long-time behavior of solutions of such
equations in a bounded domain was studied by many authors, for
instance in \cite{temam}, \cite{sell} and the references therein.

In the unbounded domain case, there also exists an extensive
literature. In 1994, E.~Feireisl \cite{feireisl1} showed that the
more challenging equation (\ref{inwav2}) admits a global attractor
 in $H^1(\R^3)\times L^2(\R^3)$ when $N=3$. For arbitrary $N$, he
obtain in \cite{feireisl2} the same result in the phase space
$H^{1}_{loc}(\R^N)\times L^{2}_{loc}(\R^N)$.

In 2001, S.V.~Zelik \cite{zelik} considered the nonautonomous case
for equation (\ref{inwav1}), in which the forcing term $g$ depends
on time. He obtained the existence of locally compact global
attractor and the upper and lower bounds for their Kolomogorov's
$\varepsilon$-entropy. Some other authors have also considered
different types of wave equations in unbounded domains
(\cite{karachalios1}, \cite{karachalios2}) in weighted spaces.

In this paper, we establish the existence of global attractors in
the usual Hilbert spaces $H^1\times L^2$, for equations
(\ref{inwav1}) and (\ref{inwav2}) in unbounded domains of $\R^N$.
To this end we cannot apply the same procedure as for bounded
domains, since the Sobolev embeddings are no longer compact. We
will  apply the "tail estimation" method, introduced for the first
order lattice systems and the reaction diffusion equations
(\cite{bates,wang,wang3}). It features an "approximation" of $R^N$
by sufficiently large bounded domains $\Omega_k$, then using the
compactness of the embeddings in $\Omega_k$ and showing the
uniform null convergence of the solutions on
$\mathbb{R}^N-\Omega_k$, we finally arrive to get the asymptotical
compactness of the semiflow.\\

Before we get into the details of our work, let us first introduce 
the basic definitions and results relevant to the general theory of
dynamical systems. We present  the
notions of semiflows and attractors along with a brief presentation
of the theory of semigroups and its relation to solving abstract
nonlinear equations in a Banach space. 

\subsection{Semiflows and Attractors}
We will use in this work the definition of semiflows as in
Temam~\cite{temam}. A stronger version can be found in
Sell\,\&\,You~\cite{sell}, where the only difference is the
continuity property.
\begin{df}
Let $(H, d)$ be a complete metric space. A family of operators
$\displaystyle\{S(t)\}_{t\geq 0}$ is called a \textbf{semiflow} on
$H$, if it satisfies the following properties:
\begin{enumerate}
\item $\displaystyle S(0)=I $  (identity in $H$), i.e. $S(0)u=u\quad\forall u\in H$,
\item $\displaystyle S(t)S(s)=S(t+s), \qquad\forall  s,t\in
\R^+$,
\item  The mapping \,\,
$\displaystyle S(t):H\rightarrow H $ \,\, is continuous for every
$t\geq 0$.
\end{enumerate}
\end{df}
We introduce now the concepts of invariant sets and attractors of
a semiflow
\begin{df}
Let $S(t)$ be a semiflow on $H$ and  $K\subset H$ . We say that
$K$ is \textbf{positively invariant} if $S(t)K\subset K$, for all
$t\geq 0$. $K$ is \textbf{invariant} if $S(t)K=K$, for all $t\geq
0$.
\end{df}
To define attractors we will need the following asymmetric
Hausdorff pseudodistance:
\begin{equation}\label{pseudometric}
\sup_{a\in A}  \mathop{\inf\vphantom{p}}_{b\in B} d(A,b)
\end{equation}
where $A,\,\,B$ are bounded sets in $H$.

We say that $A$ \textbf{attracts} $B$ if
\begin{equation}
h(S(t)B,A)\rightarrow 0,\qquad \textrm{as } t\rightarrow \infty,
\end{equation}
that is: for every $\varepsilon>0$, there exists $T\geq 0$ such
that $d(S(t)u,A)\leq \varepsilon$, for all $t\geq T$ and $u\in B$.
\begin{df}
A subset $\mathcal{A}$ of $H$ is called an \textbf{attractor} for
the semiflow $S(t)$ provided that
\begin{enumerate}
\item $\mathcal{A}$ is a compact, invariant set in $H$, and
\item there is a neighborhood $U$ of $\mathcal{A}$ in $H$ such
that $\mathcal{A}$ attracts every bounded set in $U$.
\end{enumerate}
\end{df}
An attractor $\mathcal{A}$ that attracts every bounded set in $H$
is called a \textbf{global attractor}.

The existence of global attractor is in general related to what
some authors call the ``dissipativity" of the dynamical system.
This is equivalent to the existence of absorbing sets.
\begin{df}
Let $\mathcal{B}$ be a subset of $H$ and $\mathcal{U}$ an open set
containing $\mathcal{B}$. We say that $\B$ is an \textbf{absorbing
set} in $\U$ if the orbit of any bounded set in $\U$ enters into
$\B$ after a finite time (which may depend on the set):
\begin{displaymath}
\left\{
\begin{array}{lll}
\forall \B_0\subset\U, \quad \B_0 \,\,\textrm{bounded}\\
\exists t_1(\B_0) \,\, \textrm{ such that } S(t)\B_0\subset\B,
\quad \forall t\geq t_1(\B_0).
\end{array}\right.
\end{displaymath}
We also say that $\B$ attracts the bounded sets of $\U$.
\end{df}
We have also the  related concept of asymptotical compactness.
\begin{df}
A semiflow $\{S(t)\}_{t\geq 0}$ is said to be
\textbf{asymptotically compact} on $\U$ if for every bounded
sequence $\{u_n\}$ in $\U$ and $t_n\rightarrow\infty$,
$\{S(t_n)u_n\}_{t\geq 0}$ is precompact in $H$.
\end{df}
We are now ready to present a standard result on the existence of
global attractors which can be found in \cite{hale, temam}.
\begin{theo}\label{globalat}
Let $\{S(t)\}_{t\geq 0}$ be a semiflow in $X$. If $\{S(t)\}_{t\geq
0}$ has a bounded absorbing set and is asymptotically compact in
$H$, then $\{S(t)\}_{t\geq 0}$ possesses a global attractor which
is a compact invariant set that attracts every bounded set in $H$.
\end{theo}


\subsection{Semigroups of Linear Operators}
In practice, semiflows are generated by the solutions of
differential equations. We will consider abstract nonlinear ODEs
of the form
\begin{equation}\label{nonleq}
\frac{du}{dt} +Au =F(u,t)
\end{equation}
in a Banach space $X$, where $A$ is an unbounded linear operator in
$X$ and $\displaystyle F: X\times\R\rightarrow X$ is a nonlinear
functional. In this section we will present the general existence
theory for equations such as (\ref{nonleq}). This will apply
directly to a wide range of evolutionary partial differential
equations. We will first give some basic notions on semigroup theory
which is related to solving the corresponding linear  problem
\begin{equation}\label{lineq}
\frac{du}{dt} +Au =0.
\end{equation}

In  the remainder of this section, $X$ denotes a Banach space with
norm $\|\cdot\|_X$ and $\EL(X)$ is the space of bounded linear
operators on $X$.
\begin{df}
We will say that a family of operators $\{T(t)\}_{t\geq 0}$ is a
$\mathbf{C_0}$\textbf{-semigroup of linear operators} on $X$, if
$T(t)\in \EL(X)$ for all $t\in [0,+\infty)$ and the following
hold:
\begin{itemize}
\item[(i)] $T(0)=I$ (identity in  X)
\item[(ii)] $T(t)T(s)=T(t+s),\qquad s,t\in[0,+\infty)$
\item[(iii)] $\displaystyle\lim_{t\rightarrow 0^+} T(t)x=x$,\qquad for
all  $x\in X$.
\end{itemize}
\end{df}
We see that a $C_0$-semigroup is a typical example of a semiflow
on $X$.
\begin{df}
Let $T(t)$ be a $C_0$-semigroup on $X$, its \textbf{infinitesimal
generator} is the linear operator $A$ on $X$ defined as follows
\begin{itemize}
\item The domain of $A$ is:
$$D(A)=\{x\in X\,\,:\,\,\lim_{h\rightarrow
0^+}\frac{T(h)-I}{h}x\,\,\,\,\textrm{exists in } X\}$$
\item for $x\in D(A)$ we set:
$$Ax=\lim_{h\rightarrow
0^+}\frac{T(h)-I}{h}x=\frac{d^+(T(t)x)}{dt}|_{t=0}.$$
\end{itemize}
\end{df}
Next we will give a necessary and sufficient condition for an
operator to be the infinitesimal generator of a $C_0$-semigroup in a
Hilbert space $H$. We need to introduce first some concepts. Let $H$
be a Hilbert space with inner product $\langle\cdot,\cdot\rangle$. A
linear operator $\displaystyle A:D(A)(\subset
 H)\rightarrow H$ is said to be \textbf{accretive} if
 $$\textrm{Re}\,\langle Ax,x\rangle\,\geq 0, \qquad \forall x\in D(A).$$
 If in addition we have $R(I+A)=H$ ( range of $I+A$ is equal to
 $H$) then we say that $A$ is \textbf{maximal accretive}.\\
 A $C_0$-semigroup is said to be \textbf{nonexpansive} if
 $\|T(t)\|\leq 1$ for every $t\geq 0$.
 \begin{theo}[Lumer-Phillips]\label{lumer}
Let $H$ be a Hilbert space. Then a linear operator $\displaystyle
-A:D(A)(\subset H)\rightarrow H$ is the infinitesimal generator of
a nonexpansive $C_0$-semigroup $e^{-At}$ on $H$ if and only if
both the following condtions are satisfied:
\begin{itemize}
\item[(1)] A is a closed linear operator and $D(A)$ is dense in
$H$, and
\item[(2)] A is a maximal accretive operator.
\end{itemize}
 \end{theo}
This is a classical result on semigroups and their generators. The
proof can be found in \cite{pazy, sell}. However this result applies
only to Hilbert spaces; there is a more general one on Banach
spaces, namely the \textbf{Hille-Yosida} theorem.

Now  we will briefly present  the existence theory for abstract nonlinear evolutionary equations in a Banach space $X$.
There exists a vast literature on this issue, but we will just give some basic results, see \cite{pazy, sell, temam}.

We consider the following initial value problem in the Banach space $X$:
\begin{eqnarray}\label{nonlineq2}
\left\{\begin{array}{lll} \displaystyle\frac{du}{dt}+Au=F(u)\\
u(t_0)=u_0\in X,\qquad t\geq t_0\geq 0.
\end{array}\right.
\end{eqnarray}
Assume that the nonlinearity $F$ belongs to $F\in
C_{Lip}=C_{Lip}(X,X)$, the collection of all continuous functions
$G: X\rightarrow X$ that are Lipschitz continuous on every bounded
set $B$ in $X$. We suppose also that $-A$ generates a
$C_0$-semigroup $T(t)$ on $X$.

At first, we give different notions of solution for  problem
(\ref{nonlineq2}) and then present some existence results for such
types of solutions.
\begin{df}
Let $I=[t_0,t_0+\tau)$ be an interval in $\R^+$, where $\tau>0$. A
strongly continuous mapping $u:\,I\rightarrow X$ is said to be a
\textbf{mild solution} of (\ref{nonlineq2}) in X if it solves the
following integral equation
\begin{equation}\label{inteq}
u(t)=T(t-t_0)u_0+\int_{t_0}^{t}T(t-s)F(u(s))\,ds,\qquad t\in I.
\end{equation}
If $u$ is differentiable almost everywhere in $I$ with $u_t, Au\in
L^{1}_{loc}(I,X)$, and satisfies the differential equation
\begin{equation}\label{strgsol}
\displaystyle\frac{du}{dt}+Au=^{a.e.}F(u),\,\,\, \textrm{on }
(t_0,t_0+\tau),\quad\textrm{and } u(t_0)=u_0,
\end{equation}
then $u$ is called a \textbf{strong solution} of
(\ref{nonlineq2}). If in addition, one has $u_t\in C(I,X)$ and the
differential equation in (\ref{strgsol}) is satisfied for
$t_0<t<t_0+\tau$, then $u$ is called a \textbf{classical solution}
of (\ref{nonlineq2}) on $I$.
\end{df}

We have the following result which is a particular case of
Theorems~46.1 \& 46.2 in G.~Sell \& Y.~You \cite{sell}.
\begin{theo}\label{nonlexist}
Let $-A$ generate a $C_0$-semigroup $T(t)$ on $X$ and $F\in
C_{Lip}=C_{Lip}(X,X)$. Then for every $u_0\in X$ and $t_0>0$, the
Initial Value Problem (\ref{nonlineq2}) has a unique mild solution
$u$ in $X$ on some interval $I=[t_0,t_0+\tau)$, for some $\tau>0$.

Assume that $X=H$ is a Hilbert space or a reflexive Banach space.
If $u_0\in D(A)$ or $T(t)$ is a differentiable semigroup, then the
mild solution is a strong one.
\end{theo}
\begin{rk}
The solution $u$ in Theorem~\ref{nonlexist} can be extended to a
maximum possible interval $I$. Indeed $u$ is maximally deifined if
either $\tau=+\infty$ or $\displaystyle
\lim_{t\rightarrow\tau^-}\|u(t)\|_X=+\infty$.
\end{rk}

\subsection{Some Useful  Inequalities}
We present in this section some inequalities that will be used in the remaining of this paper. The most used inequality throughout our
work is the Gronwall inequality which comes in different forms. We
present here some variants of it.
\begin{lem}[The Gronwall inequality]
Suppose that $a$ and $b$ are nonnegative constants and $u(t)$ a
nonnegative integrable function. Suppose that the following
inequality holds for $0\leq t\leq T$:
\begin{equation}
u(t)\leq a+b\int_{0}^{t}u(s)\,ds.
\end{equation}
Then for $0\leq t\leq T$, we have
\begin{equation}
u(t)\leq ae^{bt}.
\end{equation}
\end{lem}
\begin{lem}[The uniform Gronwall inequality]\label{ugronwall}
Let $g,h,y$ be nonnegative functions in $\displaystyle
L^{1}_{loc}[0,T;\R)$, where $0<T\leq\infty$. Assume that $y$ is
absolutely continuous on $(0,T)$ and that
\begin{equation}
\frac{dy}{dt}\leq gy+h \qquad \textrm{almost everywhere on }
(0,T).
\end{equation}
Then $y\in L^{\infty}_{loc}(0,T;\R)$ and one has
\begin{equation}\label{gronwineq}
y(t)\leq
y(t_0)\textrm{exp}\Big(\int_{t_0}^{t}g(s)\,ds\Big)+\int_{t_0}^{t}\textrm{exp}\Big(\int_{s}^{t}g(r)\,dr\Big)h(s)\,ds,
\end{equation}
for $0< t_0<t<T$. If in addition one has $y\in C[0,T;\R)$, then
inequality (\ref{gronwineq}) is valid at $t_0=0$.
\end{lem}
The following theorem is concerned with the Poincar\'e inequality.
\begin{theo}[Poincar\'e inequality]
Let $\Omega$ be a domain of $\R^N$ bounded only in one direction and
let $u\in H^{1}_{0}(\Omega)$. Then there is a positive constant $C$
depending only on $\Omega$ and $n$ such that
\begin{equation}
\|u\|_{L^2(\Omega)}\leq C\|\nabla u\|_{(L^2(\Omega))^N}, \qquad
\forall u\in H^{1}_{0}(\Omega).
\end{equation}
\end{theo}
\begin{rk}
The Poincar\'e inequality is usually presented for bounded domains
but the proof requires only the boundedness in one direction
$x_i$.
\end{rk}

\section{The Damped Wave Equation in $\R^N$}
We consider in this section, the nonlinear wave equation with mass
term,
\begin{eqnarray*}
u_{tt} + \lambda u_t-\triangle u+u+f(u)=g(x),  \quad t>0
\end{eqnarray*}
in $\R^N$. We shall establish first the existence and boundedness
of solutions, then we shall prove the asymptotic compactness of
the corresponding semiflow to obtain the global attractor.
\subsection{Existence of Solutions and Absorbing Set}
We start  by transforming our problem into an abstract ODE in the
space $L^2\times H^1$ and prove that the new operator is maximal
accretive. This will allow us to show the existence of solutions
and the uniform boundedness of such solutions.

We consider the system
\begin{equation}\label{wav1}
u_{tt} +\lambda u_t-\triangle u+u+f(u)=g(x), \qquad x\in
\mathbb{R}^N, \quad t>0 \\
\end{equation}
with initial conditions
\begin{equation}\label{wav2}
u(0,x)=u_0(x), \quad u_t(0,x)=u_1(x)\qquad x\in \mathbb{R}^N
\end{equation}
where $\lambda >0$,\quad $g\in L^2(\mathbb{R}^N)$, and $f\in
C^1(\mathbb{R},\,\,\mathbb{R})$ satisfies the following condition:
\begin{equation}\label{nonlincond1}
f(0)=0, \qquad f(s)s\geq \nu F(s)\geq 0, \qquad \forall
s\in\mathbb{R}
\end{equation}
where $\nu$ is a positive constant and $\displaystyle
F(s)=\int_{0}^{s}f(t)\,dt$. In addition we assume that
\begin{equation} \label{nonlincond2}
0\leq\limsup_{s \to \infty}\frac{f(s)}{s}<\infty
\end{equation}

Now set $H=L^2(\mathbb{R}^N)$, \quad $V=H^1(\mathbb{R}^N)$,\quad
and $X=V\times H$ with the usual norms and scalar products. We
define the operator $G$ in $X$ by:
\begin{eqnarray}
\begin{array}{ll}\\
D(G)=H^2(\mathbb{R}^N)\times H^1(\mathbb{R}^N)\\ \\ \\
Gw=\left(\begin{array}{ll}
\qquad \qquad\delta u-v \\ \\
-\Delta u+(\lambda-\delta)v +(\delta^2-\delta\lambda +1)u
\end{array}\right)
\end{array}
\end{eqnarray}
for $w=(u,v)\in D(G)$.\\

Then (\ref{wav1}), (\ref{wav2}) are equivalent to the initial
value problem in $X$:
\begin{eqnarray}\label{wav3}
\left\{\begin{array}{ll}
w_t+Gw=R(w), \qquad t>0, \quad w\in X \\ \\
w(0)=w_0=(u_0,\,v_0+\delta u)
\end{array}\right.
\end{eqnarray}
where
\begin{displaymath}
R(w)=\left(
\begin{array}{ll}
\qquad 0\\ \\
-f(u) +g
\end{array}\right)
\end{displaymath}

 The next result establishes the maximal accretivity of the
the operator $G$  in $X$.
\begin{lem}\label{maxac}
For a suitable $\delta$ chosen to be
$\displaystyle\delta=\frac{\lambda}{\lambda^2+4}$, the operator
$G$ defined previously is maximal accretive in $X$, and there
exists a constant $C(\delta)>0$ depending on $\delta$ such that
\begin{equation}\label{accr}
\langle Gw,\,w\rangle _X\,\,\geq C(\delta)\|w\|_{X}^{2}, \qquad
\forall w\in D(G)
\end{equation}
\end{lem}
{\bf Proof:} We first prove the positivity. Let $w=(u,\,v)\in
D(G)$, then
\begin{displaymath}
\begin{array}{lll}
\langle Gw\,,\,w\rangle _{X} & = & \langle\delta u-v,u\rangle _{V}
+\langle -\Delta u +(\lambda -\delta)v
+(\delta^2-\delta\lambda +1)u,\,v\rangle _{H}\\ \\
& = & \delta \|u\|^{2}_{V} -\int_{\mathbb{R}^N} \nabla
u\cdot\nabla v\,dx-\langle u,\,v\rangle _{H}+\int_{\mathbb{R}^N}
\nabla u\cdot\nabla v\,dx
\\ \\ &\,\, &+(\lambda
-\delta)\|v\|_{H}+(\delta^2-\delta\lambda +1)\langle u,\,v\rangle _{H}\\ \\
&=& \delta\|u\|^2_{V}+(\lambda
-\delta)\|v\|^2_{H}+(\delta^2-\delta\lambda)\langle u,\,v\rangle
_{H}
\end{array}
\end{displaymath}
Then setting
 \begin{equation}\label{sigma}
\sigma=\frac{\lambda}{\sqrt{\lambda^2+4}(\lambda+\sqrt{\lambda^2+4})},
 \end{equation}
we have
\begin{displaymath}
\begin{array}{lll}
\langle G(w)\,,\,w\rangle_X
-\sigma\big(\|u\|^{2}_{V}+\|v\|^{2}_{H}\big)-\frac{\lambda}{2}\|v\|^{2}_{H}
&\geq & (\delta-\sigma)\|u\|_{V}^{2}+(\frac{\lambda}{2}
-\delta-\sigma)\|v\|^{2}_{H}\\  & &-\delta\lambda\|u\|_V\|v\|_H\\ \\
&\geq & 2\sqrt{(\delta-\sigma)(\frac{\lambda}{2}
-\delta-\sigma)}\|u\|_V\|v\|_H\\  & &-\delta\lambda\|u\|_V\|v\|_H.
\end{array}
\end{displaymath}
We can check that $4(\delta-\sigma)(\frac{\lambda}{2}
-\delta-\sigma)=\lambda^2\delta^2$ so that
$$\langle G(w)\,,\,w\rangle_X
-\sigma\|w\|^{2}_{X}-\frac{\lambda}{2}\|v\|^2\geq 0.$$ It suffices
to take $C(\delta)=\sigma$

Now we prove that  the range of $\displaystyle G+I$ equals $X$.
Let $f=(h,g)\in X$; the question is whether there exists a
$\displaystyle w=(u,v)\in D(G)$ such that:
$$Gw+w=f\quad?$$
\begin{displaymath}
\begin{array}{ll}
\textrm{i.e.} & \left\{\begin{array}{l}
 \delta u -v+u=h       \\
-\Delta u +(\lambda -\delta)v+(\delta^2-\delta\lambda +1)u=g
\end{array}\right. \\ \\
\textrm{i.e.} & \left\{\begin{array}{l} v=(\delta+1)u-h       \\
-\Delta u +(\lambda -\delta)[(\delta+1)u-h]
+(\delta^2-\delta\lambda +1)u=g
\end{array}\right. \\ \\
\textrm{i.e.} & \left\{\begin{array}{l} v=(\delta+1)u-h        \\
 -\Delta u +(\lambda -\delta +1)u=g+(\lambda -\delta)h
\end{array}\right.
\end{array}
\end{displaymath}
Note that the operator $Au=-\Delta u$ in $L^2(\mathbb{R}^N)$ with
domain $H^2(\mathbb{R}^N)$ is a sectorial operator and there
exists $\omega\in\mathbb{R}$ such that
$\displaystyle\varrho(A)\supset\{\lambda\in\mathbb{C}\,:\,
Re\,\lambda\geq\omega\}$. So the equation
$$-\Delta u +(\lambda -\delta +1)u=g+(\lambda -\delta)h$$ has a
unique solution $u\in H^2(\mathbb{R}^N)$, thus letting
$v=(\delta+1)u-h$ and $w=(u,v)$, we get a unique $w\in D(G)$ such
that $\displaystyle Gw+w=f$. So the range of $\displaystyle G+I$
equals $X$. This, with (\ref{accr}), shows that $G$ is maximal
accretive and finishes the proof of
lemma~\ref{maxac}.\\

Lemma~\ref{maxac} together with  the Lumer-Phillips
theorem imply that $-G$ generates a nonexpansive
$C_0$-semigroup $\displaystyle e^{-Gt}$ on $X$. Furthermore since
$f$ verifies (\ref{nonlincond2}), the operator $\displaystyle R:
X\to X$ is locally Lipschitz continuous. By the standard theory of
evolutionary equations (see G.~R.~Sell \& Y.~You \cite{sell},
Theorem~46.1) this leads to the existence and uniqueness of local
solutions as stated in the next lemma.

\begin{lem}
If $g\in L^2(\mathbb{R}^N)$ and $f$ satisfies (\ref{nonlincond2}),
then for any initial data $\displaystyle w_0=(u_0,\,v_0)\in X$,
there exists a unique local solution $\displaystyle w(t)=
(u(t),\,v(t))$ of (\ref{wav3}) such that $\displaystyle w\in
C^1((-T_0,\,T_0),\,E)$ for some $T_0= T_0(w_0)
>0$.
\end{lem}

In fact we will show that the local solution $w(t)$ of
(\ref{wav3}) is bounded and  exists globally.

\begin{lem}\label{dissipativity}
Assume that (\ref{nonlincond1}) and (\ref{nonlincond2}) are
satisfied and that $g\in H$. Then any solution $w(t)$ of problem
(\ref{wav3}) satisfies
\begin{equation}\label{dissip}
\|w(t)\|_X\leq M, \qquad t\geq T_1
\end{equation}
where $M$ is a constant depending only on $\displaystyle
(\lambda,\,\, g)$ and $T_1$ depending on the data $\displaystyle
(\lambda,\,\, g,\,\, R)$ when $\|w_0\|_X\leq R$.
\end{lem}
{\bf Proof:} Let $w_0\in D(G)$ be the initial condition in
(\ref{wav3}). Taking the inner-product of (\ref{wav3}) with $w$ in
$X$ we find that
\begin{displaymath}
\begin{array}{lll}
\frac{1}{2}\frac{d}{dt}\|w\|^{2}_{X} & = & -\langle Gw\,,\,w\rangle _{X}+\langle R(w),\,w\rangle _E \\
\\& = & -\langle Gw\,,\,w\rangle _X+\langle g\,,\,v\rangle _H-\langle f\,,\,v\rangle _H\\ \\
& \leq & -C(\delta)\|w\|_{X}^{2} +\|g\|_{H}\|v\|_H -\delta\langle f(u)\,,\, u\rangle _H-\langle f(u)\,,\,u_t\rangle_H, \\ \\
\end{array}
\end{displaymath}
by (\ref{nonlincond1}) we have $$-\delta\langle f(u)\,,\, u\rangle
_H \leq -\delta\nu \int_{\mathbb{R}^N}F(u)\,dx$$ and $$-\langle
f(u)\,,\,u_t\rangle_H=-\frac{d}{dt}\int_{\mathbb{R}^N}F(u)\,dx.$$
Then using the Young inequality, it follows for any $\alpha
>0$ that
$$\frac{1}{2}\frac{d}{dt}\|w\|^{2}_{X}\,\leq\, -C(\delta)\|w\|_{X}^{2}+\frac{\alpha}{2}\|v\|^{2}_{H}+\frac{1}{2\alpha}\|g\|^{2}_{H}
-\delta\nu \int_{\mathbb{R}^N}F(u)\,dx
-\frac{d}{dt}\int_{\mathbb{R}^N}F(u)\,dx$$ which implies that
$$\frac{d}{dt}\left[\|w\|^{2}_{X}+2\int_{\mathbb{R}^N}F(u)\,dx\right]\leq 2\left (\alpha-C(\delta)\right )\|w\|^{2}_{X}-2\delta\nu \int_{\mathbb{R}^N}F(u)\,dx+\frac{1}{\alpha}\|g\|^{2}_{H}$$
Now we can choose $\alpha$ small enough so that
$\alpha-C(\delta)<0$ and taking \\$\displaystyle\mu
=\min\left\{-2(\alpha-C(\delta))\,\,,\,\,\delta\nu\right\}>0$ we
have
\begin{equation}\label{dissip1}
\frac{d}{dt}\left[\|w\|^{2}_{X}+2\int_{\mathbb{R}^N}F(u)\,dx\right]\leq
-\mu\left[\|w\|^{2}_{X}+2\int_{\mathbb{R}^N}F(u)\,dx\right]+\frac{1}{\alpha}\|g\|^{2}_{H}
\end{equation}
and then by Uniform Gronwall inequality we get
\begin{eqnarray*}
\|w\|^{2}_{X}+2\int_{\mathbb{R}^N}F(u)\,dx\leq e^{-\mu
t}\left(\|w_0\|^{2}_{X}+2\int_{\mathbb{R}^N}F(u_0)\,dx\right)+(1-e^{-\mu
t})\frac{1}{\mu\alpha}\|g\|^{2}_{H}
\end{eqnarray*}
which yields
\begin{equation}\label{dissip2}
\|w\|^{2}_{X}\leq e^{-\mu
t}\left(\|w_0\|^{2}_{X}+2\int_{\mathbb{R}^N}F(u_0)\,dx\right)+\frac{1}{\mu\alpha}\|g\|^{2}_{H}.
\end{equation}
Now by (\ref{nonlincond1}) we have $$\int_{\R^N} F(u_0)\,dx\leq
\frac{1}{\nu}\int_{\R^N}f(u_0)u_0\,dx\leq
\frac{C}{\nu}\int_{\R^N}u_{0}^{2}(x)\,dx.$$ Then we deduce from
(\ref{dissip2}) that for every $w_0\in D(G)$,
\begin{equation}\label{dissip3}
\|w\|^{2}_{X}\leq e^{-\mu
t}\left(\|w_0\|^{2}_{X}+\frac{C}{\nu}\|u_0\|_{H}^{2}\right)+\frac{1}{\mu\alpha}\|g\|^{2}_{H}.
\end{equation}
And by density of $D(G)$ in $X$ and  the continuity of the
solution of (\ref{wav3}) in $\displaystyle X\times (0\,,\,T(w_0))$
we see that (\ref{dissip2}) holds for every $w_0\in X$.\\

Now let $R>0$ and $\|w_0\|_X\leq R$, then $\|u_0\|_H\leq R$ and
\begin{equation}\label{dissip4}
\|w\|^{2}_{X}\leq e^{-\mu
t}\left(R^{2}+\frac{CR^2}{\nu}\right)+\frac{1}{\mu\alpha}\|g\|^{2}_{H}
\end{equation}
which yields
\begin{equation}\label{dissip5}
\|w\|^{2}_{X}\leq
\frac{2}{\mu\alpha}\|g\|^{2}_{H},\qquad\textrm{for } t\geq
T_1=\frac{1}{\mu}\ln\left\{
\frac{\mu\lambda(R^2+\frac{CR^2}{\nu})}{\|g\|^{2}_{H}}\right\}
\end{equation}
and (\ref{dissip}) follows with $\displaystyle
M=\frac{2}{\mu\alpha}\|g\|^{2}_{H}$ and the proof is complete.

By (\ref{dissip4}), We have also the following result.
\begin{lem}\label{bddestimate}
Let $g\in H$. Then for any given $T>0$, every solution $w$ of
(\ref{wav3}) satisfies
\begin{equation}\label{bddest}
\|w\|_X\leq L, \qquad 0\leq t\leq T
\end{equation}
where $L$ depends on $(\lambda,\delta,\|g\|_H)$, $T$ and
$\|w_0\|_X$.
\end{lem}

Lemma~\ref{dissipativity} implies that the solution $w(t)$ exists
globally, that is $T(w_0)=+\infty$, which implies that the system
(\ref{wav3}) generates a continuous semiflow
$\displaystyle\{S(t)\}_{t\geq 0}$ on $X$. Denote by $O$ the ball
\begin{equation}\label{absorb}
O=\left\{w\in X\, :\, \|w\|_X\leq M\right\}
\end{equation}
where $M$ is the constant in (\ref{dissip}). Then it follows from
(\ref{dissip}) that $O$ is an absorbing set for $S(t)$ in $X$ and
that for every bounded set $B$ in $X$ there exists a constant
$T(B)$ depending only on $(\lambda ,g)$ and $B$ such that
\begin{equation}\label{absorb1}
S(t)B\subseteq O , \quad t\geq T(B).
\end{equation}
In particular there exists a constant $T_0$ depending only on
$(\lambda ,g)$ and $O$ such that
\begin{equation}\label{absorb2}
 S(t)O\subseteq O , \quad t\geq T_0.
\end{equation}

\subsection{Global Attractor} The
existence of an absorbing set is the first step toward the existence
of a global attractor. We need now to prove the asymptotic
compactness of $S(t)$. The key idea lies in establishing uniform
estimates on ``Tail Ends" of solutions, that is, the norm of the
solutions $w(t)$ are uniformly small with respect to $t$ outside a
sufficiently large ball.
\begin{lem}\label{tail}
If (\ref{nonlincond1}) and (\ref{nonlincond2}) hold, $g\in H$ and
$w_0=(u_0,v_0)\in O$, then for every $\varepsilon >0$, there
exists positive constants $T(\varepsilon)$ and $K(\varepsilon)$
such that the solution $w(t)=(u(t),v(t))$ of problem (\ref{wav3})
satisfies
\begin{equation}\label{tail1}
\int_{|x|\geq k}\Big\{|u|^2+|\nabla u|^2+|v|^2\Big\}\,dx\leq
\varepsilon, \qquad t\geq T(\varepsilon),\quad k\geq
K(\varepsilon).
\end{equation}
\end{lem}
{\bf Proof:} Choose a smooth function $\theta$ such that $0\leq
\theta(s)\leq 1$ for $s\in \mathbb{R}^{+}$, and
$$\theta(s)=0 \quad \textrm{for } 0\leq s\leq 1; \qquad
\theta(s)=1\quad \textrm{for } s\geq 2.$$ Then there exists a
constant $C>0$ such that $|\theta'(s)|\leq C$ for $s\in
\mathbb{R}^+$.\\ Let $w(t)=(u(t),v(t))$ be the solution of problem
(\ref{wav3}) with initial condition $w_0=(u_0,v_0)\in O$ then
$v(t)=\delta u+u_t$ satisfies the equation
\begin{equation}\label{tail2}
v_t-\Delta u+ (\lambda
-\delta)v+(\delta^2-\lambda\delta+1)u=-f(u)+g
\end{equation}
taking inner product of (\ref{tail2}) with
$\theta(\frac{|x|^2}{k^2})v$ in $H$ we get
\begin{eqnarray}\label{tail3}
\begin{array}{lll}\displaystyle
\int_{\mathbb{R}^N}\theta(\frac{|x|^2}{k^2})vv_t\,dx
-\int_{\mathbb{R}^N}\Delta u\theta(\frac{|x|^2}{k^2})v\,dx +
(\lambda-\delta)\int_{\mathbb{R}^N}\theta(\frac{|x|^2}{k^2})|v|^2\,dx
\\
\\\displaystyle +(\delta^2-\lambda\delta+1)\int_{\mathbb{R}^N}\theta(\frac{|x|^2}{k^2})uv\,dx=
-\int_{\mathbb{R}^N}f(u)\theta(\frac{|x|^2}{k^2})v\,dx+\int_{\mathbb{R}^N}\theta(\frac{|x|^2}{k^2})gv\,dx
\end{array}
\end{eqnarray}
But\\
\begin{displaymath}
\begin{array}{lll}\displaystyle
-\int_{\mathbb{R}^N}\Delta u\theta(\frac{|x|^2}{k^2})v\,dx & = &
\displaystyle \int_{\mathbb{R}^N}\theta(\frac{|x|^2}{k^2})\nabla
u\cdot\nabla
v\,\,+\frac{2}{k^2}\int_{\mathbb{R}^N}\theta^\prime(\frac{|x|^2}{k^2})vx\cdot\nabla
u \\ \\
& = & \displaystyle
\int_{\mathbb{R}^N}\theta(\frac{|x|^2}{k^2})\left[\delta|\nabla
u|^2+\nabla u\cdot\nabla u_t
\right]+\frac{2}{k^2}\int_{\mathbb{R}^N}\theta^\prime(\frac{|x|^2}{k^2})vx\cdot\nabla
u \\ \\
&=& \displaystyle
\frac{1}{2}\frac{d}{dt}\int_{\mathbb{R}^N}\theta(\frac{|x|^2}{k^2})|\nabla
u|^2+\delta\int_{\mathbb{R}^N}\theta(\frac{|x|^2}{k^2})|\nabla
u|^2
\\ \\ & &\displaystyle +\frac{2}{k^2}\int_{\mathbb{R}^N}\theta^\prime(\frac{|x|^2}{k^2})vx\cdot\nabla
u, \\
\end{array}
\end{displaymath}
and
\begin{displaymath}
\begin{array}{lll}\displaystyle
(\delta^2-\lambda\delta+1)\int_{\mathbb{R}^N}\theta(\frac{|x|^2}{k^2})uv\,dx
=(\delta^2-\lambda\delta+1)\int_{\mathbb{R}^N}\theta(\frac{|x|^2}{k^2})(\delta
|u|^2+uu_t) \\ \\ = \displaystyle
\frac{1}{2}(\delta^2-\lambda\delta+1)\frac{d}{dt}\int_{\mathbb{R}^N}\theta(\frac{|x|^2}{k^2})
|u|^2+\delta(\delta^2-\lambda\delta+1)\int_{\mathbb{R}^N}\theta(\frac{|x|^2}{k^2})
|u|^2.
\end{array}
\end{displaymath}\\
Then (\ref{tail3}) becomes\\
\begin{eqnarray}\label{tail4}
\begin{array}{lll}\displaystyle
\frac{1}{2}\frac{d}{dt}\int_{\mathbb{R}^N}\theta(\frac{|x|^2}{k^2})\left[(\delta^2-\lambda\delta+1)|u|^2+|\nabla
u|^2+|v|^2\right]\\
\\ \displaystyle+\delta\int_{\mathbb{R}^N}\theta(\frac{|x|^2}{k^2})\left[(\delta^2-\lambda\delta+1)|u|^2+|\nabla
u|^2+|v|^2\right]+
(\lambda-2\delta)\int_{\mathbb{R}^N}\theta(\frac{|x|^2}{k^2})|v|^2\\
\\\displaystyle=-\int_{\mathbb{R}^N}\theta(\frac{|x|^2}{k^2})f(u)(\delta u+u_t)
+\int_{\mathbb{R}^N}\theta(\frac{|x|^2}{k^2})gv\,dx
-\frac{2}{k^2}\int_{\mathbb{R}^N}\theta^\prime(\frac{|x|^2}{k^2})vx\cdot\nabla
u,
\end{array}
\end{eqnarray}\\
and since\\
$$\int_{\mathbb{R}^N}\theta(\frac{|x|^2}{k^2})f(u)(\delta
u+u_t)\geq
\frac{d}{dt}\int_{\mathbb{R}^N}\theta(\frac{|x|^2}{k^2})F(u)+\delta\nu\int_{\mathbb{R}^N}\theta(\frac{|x|^2}{k^2})F(u),$$\\
we deduce that\\
\begin{eqnarray}\label{tail5}
\begin{array}{lll}\displaystyle
\frac{d}{dt}\int_{\mathbb{R}^N}\theta(\frac{|x|^2}{k^2})\left[(\delta^2-\lambda\delta+1)|u|^2+|\nabla
u|^2+|v|^2+2F(u)\right]\\
\\\displaystyle+\delta\alpha\int_{\mathbb{R}^N}\theta(\frac{|x|^2}{k^2})\left[(\delta^2-\lambda\delta+1)|u|^2+|\nabla
u|^2+|v|^2+2 F(u)\right]\\ \\ \displaystyle \leq
-(\lambda-2\delta)\int_{\mathbb{R}^N}\theta(\frac{|x|^2}{k^2})|v|^2
+\int_{\mathbb{R}^N}\theta(\frac{|x|^2}{k^2})gv\,dx
-\frac{2}{k^2}\int_{\mathbb{R}^N}\theta^\prime(\frac{|x|^2}{k^2})vx\cdot\nabla
u,
\end{array}
\end{eqnarray}
where $\displaystyle \alpha=\min\left\{1\,\,,\,\,\nu\right\}$.
Now, there exists a constant $K(\varepsilon)>0$ such that for
$k\geq
K$, we have\\
$$\displaystyle-(\lambda-2\delta)\int_{\mathbb{R}^N}\theta(\frac{|x|^2}{k^2})|v|^2
+\int_{\mathbb{R}^N}\theta(\frac{|x|^2}{k^2})gv\,dx
-\frac{2}{k^2}\int_{\mathbb{R}^N}\theta^\prime(\frac{|x|^2}{k^2})vx\cdot\nabla
u\leq \frac{\varepsilon}{2},$$\\ which implies by Uniform Gronwall
inequality
that\\
\begin{eqnarray*}\label{tail6}
&
&\int_{\mathbb{R}^N}\theta(\frac{|x|^2}{k^2})\left[(\delta^2-\lambda\delta+1)|u|^2+|\nabla
u|^2+|v|^2+2F(u)\right]\nonumber\\
\nonumber\\& &\displaystyle\leq e^{-\delta\alpha
t}\int_{\mathbb{R}^N}\theta(\frac{|x|^2}{k^2})\left[(\delta^2-\lambda\delta+1)|u_0|^2+|\nabla
u_0|^2+|v_0|^2+2F(u_0)\right]+\varepsilon\frac{1-e^{-\delta\alpha}}{2\delta\alpha}
\nonumber.\\
\end{eqnarray*}
Now since $w_0\in O$,  there exist a constant $M>0$, uniformly
chosen for $w_0\in O$, such that
\\$$\int_{\mathbb{R}^N}\theta(\frac{|x|^2}{k^2})\left[(\delta^2-\lambda\delta+1)|u_0|^2+|\nabla
u_0|^2+|v_0|^2+2F(u_0)\right]\leq M.$$ \\Then we get for $k\geq
K(\varepsilon)$,\\
$$\int_{\mathbb{R}^N}\theta(\frac{|x|^2}{k^2})\left[(\delta^2-\lambda\delta+1)|u|^2+|\nabla
u|^2+|v|^2+2F(u)\right]\leq Me^{-\delta\alpha t}+
\varepsilon\frac{1-e^{-\delta\alpha}}{2\delta\alpha}.$$\\
Choosing $\displaystyle
T(\varepsilon)=\frac{1}{\delta}\ln\Big(\frac{2M\delta\alpha}{2\varepsilon\delta\alpha-\varepsilon}\Big)$,
we deduce that\\
$$\int_{\mathbb{R}^N}\theta(\frac{|x|^2}{k^2})\left[(\delta^2-\lambda\delta+1)|u|^2+|\nabla
u|^2+|v|^2\right]\leq\varepsilon\qquad\textrm{for }t\geq
T(\varepsilon), \,\,\, k\geq K(\varepsilon)$$\\ which yields
(\ref{tail1}) since $0<\delta^2-\lambda\delta+1<1$ for the
particular choice of $\delta$, and the proof is complete.\\

By multiplying equation (\ref{tail2}) with $v$ and integrating we
deduce the following energy equation
\begin{equation}\label{energy}
\frac{d}{dt}E(w(t))+2\delta E(w(t))=G(w(t))\qquad \forall t>0,
\end{equation}
where $E(w)$ is the quasi-energy functional,
\begin{equation}\label{energy1}
E(w)=(\delta^2-\lambda\delta+1)\|u\|^{2}_{H}+\|\nabla
u\|^{2}_{L^2(\mathbb{R}^N)}+\|v\|^{2}_{H},
\end{equation}
and
\begin{equation}\label{energy2}
G(w)=-2(\lambda-2\delta)\|v\|_{H}^{2}+2\int_{\mathbb{R}^N}gv\,dx-2\int_{\mathbb{R}^N}f(u)v\,dx.
\end{equation}
This energy functional $E$ will be used later as an equivalent norm,
more suitable in proving the asymptotical compactness.\\

The following lemma will be also  useful in proving the asymptotical
compactness.
\begin{lem}
Let $\displaystyle w_n=(u_n,v_n)\longrightarrow w_0=(u_0,v_0)$
weakly in $X$, then for every $T>0$ we have
\begin{equation}\label{weakconv1}
S(t)w_n\longrightarrow S(t)w_0\qquad\textrm{weakly in } L^2(0,T;X)
\end{equation} and
\begin{equation}\label{weakconv2}
S(t)w_n\longrightarrow S(t)w_0\qquad\textrm{weakly in } X,
\quad\textrm{for } 0\leq t\leq T.
\end{equation}
\end{lem}
{\bf Proof:} Since $\{w_n\}_n$ converges weakly in $X$, then  it
is bounded in $X$ so that, by lemma~\ref{bddestimate}
$\{S(t)w_n\}_n$ is bounded in $L^\infty (0,T;X)$. This, with
(\ref{wav3}),  implies
 that
\begin{equation}\label{weakconv3}\frac{\partial}{\partial
t}S(t)v_n \quad \textrm{is bounded in }
L^\infty(0,T;H^{-1}(\mathbb{R}^N))\end{equation} and
\begin{equation}\label{weakconv3}S(t)v_n \quad \textrm{is bounded in } L^\infty(0,T;L^{2
}(\mathbb{R}^N)).\end{equation} We infer that there exists a
subsequence $\{w_{n_{j}}\}_j$ and $w_\infty=(u_\infty,
v_\infty)\in L^\infty(0,T;X)$ such that
\begin{equation}\label{weakconv4}
S(t)w_{n_j}\longrightarrow w_\infty\qquad\textrm{weakly in }
L^2(0,T;X),
\end{equation}
\begin{equation}\label{weakconv5}\frac{\partial}{\partial
t}S(t)v_{n_j}\longrightarrow  \frac{\partial}{\partial t}v_\infty
\qquad\textrm{weakly in }
L^\infty(0,T;H^{-1}(\mathbb{R}^N))\end{equation} and
\begin{equation}\label{weakconv6}\frac{\partial}{\partial
t}S(t)u_{n_j}\longrightarrow  \frac{\partial}{\partial t}u_\infty
\qquad\textrm{weakly in } L^2(0,T;H).
\end{equation} We can show
that $w_\infty$ is a solution of (\ref{wav3}) with
$w_\infty(0)=w_0$. Indeed, we have by the mild solution formula,
\begin{equation}
S(t)w_{n_j}=e^{-Gt}w_{n_j}+\int_{0}^{t}e^{-G(t-s)}R(S(s)w_{n_j})\,ds.
\end{equation}
And, since $\displaystyle w_{n_j}\rightarrow w_0$ weakly in $X$, we
deduce from (\ref{weakconv4}) that
\begin{equation}
e^{-Gt}w_{n_j}+\int_{0}^{t}e^{-G(t-s)}R(S(s)w_{n_j})\,ds\longrightarrow
e^{-Gt}w_{0}+\int_{0}^{t}e^{-G(t-s)}R(w_\infty(s))\,ds,
\end{equation}
weakly in $X$. which implies, by the uniqueness of weak limit that
\begin{equation}
w_\infty(t)=e^{-Gt}w_{0}+\int_{0}^{t}e^{-G(t-s)}R(w_\infty(s))\,ds.
\end{equation}
That is $w_\infty$ is a solution of (\ref{wav3}) and by the
uniqueness of solutions we have $w_\infty (t)=S(t)w_0$. This shows
that any subsequence of $S(t)w_n$ has a weakly convergent
subsequence in $L^2(0,T;X)$, therefore we conclude
(\ref{weakconv1}). A similar argument yields (\ref{weakconv2}).

Similar to (\ref{weakconv1}) we also have that if
$w_n\longrightarrow w$ weakly in $X$, then for $0\leq s\leq T$,
\begin{equation}\label{weakconv7}
S(t)w_n\longrightarrow S(t)w_0\qquad\textrm{weakly in } L^2(s,T;X)
\end{equation}

We state here another useful lemma.
\begin{lem}\label{weakconv8}
Let $\Omega$ be a bounded domain in $\mathbb{R}^N$. Suppose
$u_n\longrightarrow u$  in $L^2(\Omega)$ and $v_n\longrightarrow
v$ weakly in $L^2(\Omega)$ , then
$\displaystyle\int_{\Omega}f(u_n)v_n\,dx\longrightarrow
\int_{\Omega}f(u)v\,dx$ in $\mathbb{R}$ (up to a subsequence
extraction).
\end{lem}
{\bf Proof:} By (\ref{nonlincond2}) we can show, up to a
subsequence, that  $f(u_n)\longrightarrow f(u)$ in $L^2(\Omega)$.
Now define the linear functionals  $I_n$ and $I$ on $L^2(\Omega)$
by
$$ I_n(v)=\int_{\Omega}f(u_n)v\,dx, \quad
I(v)=\int_{\Omega}f(u)v\,dx. $$ Then $I_n\longrightarrow I$ in
$L^2(\Omega)^*$ (the dual space of $L^2(\Omega)$. Indeed
\begin{eqnarray*}
|I_n(v)-I(v)|& \leq & \int_{\Omega}|f(u_n)-f(u)||v|\,dx\\ \\
& \leq & \|f(u_n)-f(u)\|_{L^2}\|v\|_{L^2}.
\end{eqnarray*}
which implies that
$$\|I_n-I\|_{L^2}\leq \|f(u_n)-f(u)\|_{L^2}\longrightarrow 0 \quad \textrm{as } n\longrightarrow \infty.$$
So $I_n\rightarrow I$ in $L^2(\Omega)^*$ and $v_n\rightarrow v$
weakly in $L^2(\Omega)$, then it follows, by a classical result in
functional analysis that
$I_n(v_n) \longrightarrow I(v)$, which proves the lemma. \\

We are now ready to prove the asymptotic compactness of the
semiflow $S(t)$.
\begin{theo}\label{assymp}
The semiflow $S(t)$ generated by the system (\ref{wav3}) is
asymptotically compact in $X$, that is if $\{w_n\}_n$ is a bounded
sequence in $X$ and $t_n\longrightarrow +\infty$, then
$\{S(t)w_n\}_{n\geq 1}$ is precompact in $X$.
\end{theo}
{\bf Proof:} Let $w_n$ be a bounded sequence in $X$ with
$\|w_n\|_X\leq R$ and $t_n\longrightarrow +\infty$ then by
(\ref{absorb1}) there exists a constant $T(R)>0$ depending only on
$R>0$ such that
\begin{equation}\label{assymp1} S(t)w_n\in O, \qquad \forall n\geq 1,
\quad\forall t\geq T(R).\end{equation} Since $t_n\longrightarrow
+\infty$, there exists $N_1(R)$ such that $n\geq N_1$ implies
$t_n\geq T(R)$ so that
\begin{equation}\label{assymp2}
S(t_n)w_n\in O, \qquad \forall n\geq N_1(R).
\end{equation}
Then there exists $w\in X$ such that, up to a subsequence
\begin{equation}\label{assymp3}
S(t_n)w_n\longrightarrow w\qquad \textrm{weakly in } X.
\end{equation}
Now for every $T>0$ there exists $N_2(R,T)$ such that for $n\geq
N_2(R,T)$ we have $t_n-T\geq T(R)$ so that
\begin{equation}\label{assymp4}
S(t_n-T)w_n\in O \qquad \forall n\geq N_2(R,T).
\end{equation}
Thus there is a $w_T\in O$ such that
\begin{equation}\label{assymp5}
S(t_n-T)w_n\longrightarrow w_T\qquad\textrm{ weakly in } X,
\end{equation}
and by the weak continuity (\ref{weakconv2}) we must have
$w=S(T)w_T$ which implies that
\begin{equation}\label{liminf}
\liminf_{n\rightarrow \infty} \|S(t_n)w_n\|_X\geq \|w\|_X.
\end{equation}
So we only need to prove that
\begin{equation}\label{limsup}
\limsup_{n\rightarrow \infty} \|S(t_n)w_n\|_X\leq \|w\|_X.
\end{equation}\\

By the energy equation (\ref{energy}), it follows that any
solution $w(t)=S(t)w$ of (\ref{wav3}) satisfies
\begin{equation}\label{energy3}
E(S(t)w)=e^{-2\delta(t-s)}E(S(s)w)+\int^{t}_{s}e^{-2\delta(t-r)}G(S(r)w)\,dr\,,\quad
t\geq s\geq 0.
\end{equation}
where $E$ and $G$ are given by (\ref{energy1}) and
(\ref{energy2}), respectively.

In the following,  $T_0$ is the constant in (\ref{absorb2}), and
for $\varepsilon
>0$, $T(\varepsilon)$ is the constant in (\ref{tail1}). Let
$T_0(\varepsilon)$ be a fixed constant such that $\displaystyle
T_0(\varepsilon)\geq \max\{T(\varepsilon)\,,\,T_0\}$. Taking
$T\geq T_0(\varepsilon)$, and applying (\ref{energy3}) to the
solution $S(t)(S(t_n-T)w_n)$ with $s=T_0$ and $t=T$, then we get,
for $n\geq N_2(R,T)$,
\begin{eqnarray}\label{energy4}
E(S(t_n)w_n)& = & E(S(T)(S(t_n-T)w_n))\nonumber\\
& = & e^{-2\delta(T-T_0)}E(S(T_0)(S(t_n-T)w_n))\nonumber \\
&+&\int^{T}_{T_0}e^{-2\delta(T-r)}G(S(r)(S(t_n-T)w_n))\,dr.
\end{eqnarray}
Since $T_0\geq T_0$ we have $S(T_0)(S(t_n-T)w_n)\in O$ for $n\geq
N_2(R,T)$, therefore by the definition of $E$ we find that
\begin{equation}\label{assymp6}
e^{-2\delta(T-T_0)}E(S(T_0)(S(t_n-T)w_n))\leq C
e^{-2\delta(T-T_0)},\quad   \forall n\geq N_2(R,T).
\end{equation}
On the other hand, we have
\begin{eqnarray}\label{asymp7}
\int^{T}_{T_0}e^{-2\delta(T-r)}G(S(r)(S(t_n-T)w_n))\,dr\nonumber\\
=-2(\lambda-\delta)\int^{T}_{T_0}e^{-2\delta(T-r)}\|S(r)S(t_n-T)v_n\|^2\,dr\nonumber\\
+2\int^{T}_{T_0}e^{-2\delta(T-r)}\int_{\mathbb{R}^N}
gS(r)S(t_n-T)v_n\,dxdr\nonumber\\
-2\int^{T}_{T_0}e^{-2\delta(T-r)}\int_{\mathbb{R}^N}f(S(r)S(t_n-T)u_n)S(r)S(t_n-T)v_n\,dxdr
\end{eqnarray}
Let's handle the first and last term of (\ref{asymp7}). Since we
have,
$$e^{-2\delta(T-r)}S(r)S(t_n-T)v_n\longrightarrow
e^{-2\delta(T-r)}S(r)v \textrm{ weakly in } L^2(T_0,T;H),$$ it
follows that:
\begin{eqnarray*}
\liminf_{n\rightarrow\infty}\|e^{-2\delta(T-r)}S(r)S(t_n-T)v_n\|_{L^2(T_0,T;H)}\geq
\|e^{-2\delta(T-r)}S(r)v\|_{L^2(T_0,T;H)},
\end{eqnarray*}
which implies that
\begin{eqnarray}\label{assymp8}
\limsup_{n\rightarrow\infty}-2(\lambda-\delta)\|e^{-2\delta(T-r)}S(r)S(t_n-T)v_n\|_{L^2(T_0,T;H)}\nonumber\\\leq
-2(\lambda-\delta)\|e^{-2\delta(T-r)}S(r)v\|_{L^2(T_0,T;H)}.
\end{eqnarray}
Also by (\ref{assymp5}) and (\ref{weakconv7}) we have\\
\begin{eqnarray}\label{assymp9}
\int^{T}_{T_0}e^{-2\delta(T-r)}\int_{\mathbb{R}^N}
gS(r)S(t_n-T)v_n\,dxdr\longrightarrow
\int^{T}_{T_0}e^{-2\delta(T-r)}\int_{\mathbb{R}^N} gS(r)v_T\,dxdr
\end{eqnarray}\\
Now let's handle the nonlinear term of (\ref{asymp7}). We have\\
\begin{eqnarray}\label{assymp10}
-2\int^{T}_{T_0}e^{-2\delta(T-r)}\int_{\mathbb{R}^N}f(S(r)S(t_n-T)u_n)S(r)S(t_n-T)v_n\,dxdr\nonumber\\
=-2\int^{T}_{T_0}e^{-2\delta(T-r)}\int_{|x|\geq
k}f(S(r)S(t_n-T)u_n)S(r)S(t_n-T)v_n\,dxdr\nonumber\\
-2\int^{T}_{T_0}e^{-2\delta(T-r)}\int_{|x|\leq
k}f(S(r)S(t_n-T)u_n)S(r)S(t_n-T)v_n\,dxdr.
\end{eqnarray}\\
Handling the first term on the right-hand side of (\ref{assymp10})
gives\\
\begin{eqnarray}\label{assymp11}
\left|2\int^{T}_{T_0}e^{-2\delta(T-r)}\int_{|x|\geq
k}f(S(r)S(t_n-T)u_n)S(r)S(t_n-T)v_n\,dxdr\right|\nonumber\\
\leq C\int^{T}_{T_0}e^{-2\delta(T-r)}\int_{|x|\geq
k}|S(r)S(t_n-T)u_n||S(r)S(t_n-T)v_n|\nonumber\\
\leq C\int^{T}_{T_0}e^{-2\delta(T-r)}\left(\int_{|x|\geq
k}|S(r)S(t_n-T)u_n|^2\right)^{\frac{1}{2}}\left(\int_{|x|\geq
k}|S(r)S(t_n-T)v_n|^2\right)^{\frac{1}{2}}\nonumber\\
\leq \varepsilon^2C\int^{T}_{T_0}e^{-2\delta(T-r)}\,dr\leq
\frac{\varepsilon^2C}{2\delta},\qquad\qquad n\geq N_2(R,T).\qquad
\end{eqnarray}\\
We treat now the second term on the right-hand side of
(\ref{assymp10}). We want to prove that as $n\rightarrow +\infty$,
\begin{eqnarray}\label{assymp12}
\int^{T}_{T_0}e^{-2\delta(T-r)}\int_{|x|\leq
k}f(S(r)S(t_n-T)u_n)S(r)S(t_n-T)v_n\,dxdr \nonumber\\
\longrightarrow \int^{T}_{T_0}e^{-2\delta(T-r)}\int_{|x|\leq
k}f(S(r)u_T)S(r)v_T\,dxdr\qquad\qquad
\end{eqnarray}
Set $\displaystyle\Omega_k=\{x\in\mathbb{R}^N\,:\,|x|\leq k\}$ and
let $r\in [T_0,T]$. Then we have
$$S(r)S(t_n-T)w_n\longrightarrow S(r)w_T,\qquad \textrm{weakly in
} X.$$ By the compactness of the Sobolev embedding
$H^1(\Omega_k)\subset L^2(\Omega_k)$, we infer that
\begin{equation}\label{assymp13}
S(r)S(t_n-T)u_n\longrightarrow S(r)u_T,\qquad \textrm{strongly in
} L^2(\Omega_k)
\end{equation}
and
\begin{equation}\label{assymp14}
S(r)S(t_n-T)v_n\longrightarrow S(r)v_T,\qquad \textrm{weakly in }
L^2(\Omega_k)
\end{equation}
then (\ref{assymp12}) follows from lemma~\ref{weakconv8}.\\

By (\ref{assymp10}), (\ref{assymp11}) and (\ref{assymp12}) we find
that for $k\geq K(\varepsilon)$,
\begin{eqnarray*}
\limsup_{n\rightarrow
\infty}-2\int^{T}_{T_0}e^{-2\delta(T-r)}\int_{\mathbb{R}^N}f(S(r)S(t_n-T)u_n)S(r)S(t_n-T)v_n\,dxdr\nonumber
\\
\leq \varepsilon C-2\int^{T}_{T_0}e^{-2\delta(T-r)}\int_{|x|\leq
k}f(S(r)u_T)S(r)v_T\,dxdr.
\end{eqnarray*}
Letting $k\longrightarrow \infty$ we obtain
\begin{eqnarray}\label{assymp15}
\limsup_{n\rightarrow
\infty}-2\int^{T}_{T_0}e^{-2\delta(T-r)}\int_{\mathbb{R}^N}f(S(r)S(t_n-T)u_n)S(r)S(t_n-T)v_n\,dxdr\nonumber
\\
\leq \varepsilon
C-2\int^{T}_{T_0}e^{-2\delta(T-r)}\int_{\mathbb{R}^N}f(S(r)u_T)S(r)v_T\,dxdr.
\end{eqnarray}
By (\ref{asymp7}), (\ref{assymp8}), (\ref{assymp9}) and
(\ref{assymp15}), we finally obtain
\begin{eqnarray*}
\limsup_{n\rightarrow\infty}\int^{T}_{T_0}e^{-2\delta(T-r)}G(S(r)(S(t_n-T)w_n))\,dr\nonumber\\
\leq -2(\lambda-\delta)\int^{T}_{T_0}e^{-2\delta(T-r)}\|S(r)v_T\|^2\,dr\nonumber\\
+2\int^{T}_{T_0}e^{-2\delta(T-r)}\int_{\mathbb{R}^N}
gS(r)v_T\,dxdr\nonumber\\
-2\int^{T}_{T_0}e^{-2\delta(T-r)}\int_{\mathbb{R}^N}f(S(r)u_T)S(r)v_T\,dxdr
+\varepsilon C,
\end{eqnarray*}
that is
\begin{eqnarray}\label{assymp16}
\limsup_{n\rightarrow\infty}\int^{T}_{T_0}e^{-2\delta(T-r)}G(S(r)(S(t_n-T)w_n))\,dr
\nonumber \\\leq\int^{T}_{T_0}e^{-2\delta(T-r)}G(S(r)w_T)\,dr
+\varepsilon C.
\end{eqnarray}
Taking limit of (\ref{energy4}), (\ref{assymp6}) and
(\ref{assymp16}) we get, as $n\rightarrow \infty$,
\begin{eqnarray}\label{assymp17}
\limsup_{n\rightarrow \infty} E(S(t_n)w_n) \leq Ce^{-2\delta
(T-T_0)} +\int^{T}_{T_0}e^{-2\delta(T-r)}G(S(r)w_T)\,dr
+\varepsilon C.
\end{eqnarray}
On the other hand, since $w=S(T)w_T$, by (\ref{energy3}) we also
have that
\begin{eqnarray}\label{assymp18}
E(w)=E(S(T)w_T) = e^{-2\delta (T-T_0)}E(S(T_0)w_T)
+\int^{T}_{T_0}e^{-2\delta(T-r)}G(S(r)w_T)\,dr.
\end{eqnarray}
Hence it follows from (\ref{assymp17})-(\ref{assymp18}) that
\begin{eqnarray}\label{assymp19}
\limsup_{n\rightarrow \infty} E(S(t_n)w_n) \leq E(w) +Ce^{-2\delta
(T-T_0)}+\varepsilon C -e^{-2\delta (T-T_0)}E(S(T_0)w_T).
\end{eqnarray}
Now since $w_T\in O$ and $T_0\geq T(O)$ we find that
$$|e^{-2\delta (T-T_0)}E(S(T_0)w_T)|\leq Ce^{-2\delta
(T-T_0)}.$$ Then from (\ref{assymp19}) we have
\begin{eqnarray}\label{assymp20}
\limsup_{n\rightarrow \infty} E(S(t_n)w_n) \leq E(w) +Ce^{-2\delta
(T-T_0)}+\varepsilon C.
\end{eqnarray}
Now taking limit of (\ref{assymp20}) as $T\rightarrow \infty$ and
then letting $\varepsilon\rightarrow 0$, we obtain
\begin{eqnarray*}
\limsup_{n\rightarrow \infty} E(S(t_n)w_n) \leq E(w),
\end{eqnarray*}
that is
\begin{eqnarray}\label{assymp21}
\limsup_{n\rightarrow\infty}\,\,(\delta^2-\lambda\delta+1)\|S(t_n)u_n\|^{2}_{H}+\|\nabla
S(t_n)u_n\|^{2}_{L^2(\mathbb{R}^N)}+\|S(t_n)v_n\|^{2}_{H}
\nonumber \\ \leq (\delta^2-\lambda\delta+1)\|u\|^{2}_{H}+\|\nabla
u\|^{2}_{L^2(\mathbb{R}^N)}+\|v\|^{2}_{H}.
\end{eqnarray}
Noting that $\displaystyle
E(w)=(\delta^2-\lambda\delta+1)\|u\|^{2}_{H}+\|\nabla
u\|^{2}_{L^2(\mathbb{R}^N)}+\|v\|^{2}_{H}$ is equivalent to the
norm of $X$, we can assume without loss of generality that the
norm of $X$ is defined by it. Then we have $$\limsup_{n\rightarrow
\infty} \|S(t_n)w_n\|_X\leq \|w\|_X$$ as desired in
(\ref{limsup}). Therefore we get the strong convergence of
$S(t_n)w_n$ to $w$ in $X$. The proof is complete.

Now we state our main result obtained in  this section.
\begin{theo}\label{attractor}
Assume that $f$ satisfies (\ref{nonlincond1}), (\ref{nonlincond2})
and $g\in L^2(\R^N)$. Then, problem (\ref{wav3}) possesses a
global attractor in $X=H^1(\R^N)\times L^2(\R^N)$ which is a
compact invariant subset that attracts every bounded set of $X$
with respect to the norm topology.
\end{theo}
{\bf Proof:} Since we have established the existence of an
absorbing set in (\ref{absorb}) and the asymptotic compactness of
the semiflow $S(t)$ in $X$ in Theorem~\ref{assymp}, the conclusion
follows from Theorem~\ref{globalat}.

\section{The Wave Equation Without Mass Term}
In this section we will study the existence of global attractor
for the wave equation without mass term,
\begin{eqnarray}\label{wav2.1}\left\{
\begin{array}{lll}
u_{tt} + \lambda u_t-\triangle u+f(u)=0, \qquad x\in \Omega, \quad
t>0,\\
u|_{\partial\Omega}=0,\\
u(0,x)=u_0(x),\quad u_t(0,x)=u_1(x)
\end{array}\right.
\end{eqnarray}
where $\Omega$ is a domain of $\R^N$ bounded only in one direction,
with smooth boundary. The case $\Omega=\R^N$, for this equation is
still an open problem due to some difficulties in getting an
inequality such as (\ref{accr}) for the operator $G$ in $H^1$ norm.
In our case we will use an equivalent norm (provided by the
Poincar\'e inequality) for which the desired estimate works. We
assume the same condtions (\ref{nonlincond1}) and
(\ref{nonlincond2}) for the nonlinear function $f$.

We will work in the phase space $X=V\times H$ where
$V=H^{1}_{0}(\Omega),\,\,\,H=L^2(\Omega)$. $H$ is endowed with the
 norm and inner product for $L^2$ and $V$ is endowed with the
inner product and norm defined as follows,
\begin{equation}\label{h1norm}
(u,v)_V=\int_{\Omega}\nabla u\cdot\nabla v\,dx ,\quad u,v\in
V\quad\textrm{ and } \|u\|_V=\|\nabla u\|_{H^n},\,\,\,u\in V.
\end{equation}
Now define the following bilinear operator in $V$:
\begin{equation}
(u,v)_{1}=\int_{\Omega}uv\,dx+\int_{\Omega}\nabla u\cdot\nabla
v\,dx ,\quad u,v\in V,
\end{equation}
which is also an inner product in $V$ with induced norm
$\displaystyle\|u\|_{1}=\Big[\|u\|_{H}^{2}+\|\nabla
u\|_{H^n}^{2}\Big]^{\frac{1}{2}}$. By the Poincar\'e inequality
$\|\cdot\|_V$ and $\|\cdot\|_1$ are equivalent norms in $V$.That
is, there are positive constants $C_1$ and $C_2$ such that
\begin{equation}\label{equivh1norm}
C_1\|u\|_V\leq \|u\|_1\leq C_2\|u\|_V,\qquad \forall u\in V.
\end{equation}

Let's make a transformation to write the equation (\ref{wav2.1})
as a first order abstract ODE.\\ Choose
$\displaystyle\delta=\frac{\lambda}{\lambda^2+4}$ and set
$v=\delta u+u_t$, $w=\left(\begin{array}{ll} u\\
v\end{array}\right)$. Then, problem (\ref{wav2.1}) is equivalent
to
\begin{eqnarray}\label{wav2.2}
\left\{\begin{array}{ll} w_t+Gw=R(w), \qquad t>0, \quad w\in X \\
 w(0)=w_0=(u_0,\,u_1+\delta u_0)
\end{array}\right.
\end{eqnarray}
where
\begin{displaymath}
R(w)=\left(
\begin{array}{ll}
\qquad 0\\  -f(u)
\end{array}\right)
\end{displaymath}
and
\begin{displaymath}
\begin{array}{ll}\\
Gw=\left(\begin{array}{lll} \qquad \qquad\delta u-v \\  -\Delta
u+(\lambda-\delta)v +(\delta^2-\delta\lambda )u\end{array}\right)
\end{array}
\end{displaymath}\\
for $w=\left(\begin{array}{ll} u\\
v\end{array}\right)\in D(G)=\big(H^2(\Omega)\cap H^{1}_{0}(\Omega)\big)\times H^1(\Omega)$.\\

As in Lemma~\ref{maxac} we show  the positivity of the operator
$G$ with a similar estimate.
\begin{lem}\label{maxac2}
For $\displaystyle\delta=\frac{\lambda}{\lambda^2+4}$, the operator
$G$ is maximal accretive in $X$ and verifies the following
\begin{equation}\label{acrr2}
(G(w)\,,\,w)_X\geq \sigma\|w\|_{X}^{2}+\frac{\lambda}{2}\|v\|^{2}_{H},\qquad\forall \,\,w=\left(\begin{array}{ll} u\\
v\end{array}\right)\in X,
\end{equation}
where \begin{equation}
\sigma=\frac{\lambda}{\sqrt{\lambda^2+4}(\lambda+\sqrt{\lambda^2+4})}.
\end{equation}
\end{lem}
{\bf Proof:} Let $w=\left(\begin{array}{ll} u\\
v\end{array}\right)\in X$ then we have:
\begin{displaymath}
\begin{array}{lll}
(G(w)\,,\,w)_X&=&(\delta
u-v\,,\,u)_V+(-\Delta u+(\lambda-\delta)(v-\delta u)\,,\,v)_H\\ \\
&=&\delta\|u\|^{2}_{V}-(\nabla u\,,\,\nabla v)_{H^n}+(-\Delta u\,,\,v)_H+(\lambda-\delta)\|v\|^{2}_{H}\\
\\&
&-\delta(\lambda-\delta)(u\,,\,v)_H\\ \\
&=&\delta\|u\|^{2}_{V}+(\lambda-\delta)\|v\|^{2}_{H}-\delta(\lambda-\delta)(u\,,\,v)_H\\ \\
&\geq & \delta\|u\|_{V}^{2}+(\lambda
-\delta)\|v\|^{2}_{H}-\delta\lambda\|u\|_H\|v\|_H.\\ \\
\end{array}
\end{displaymath}
Then setting
$\displaystyle\sigma=\frac{\lambda}{\sqrt{\lambda^2+4}(\lambda+\sqrt{\lambda^2+4})}$
as in (\ref{sigma}), we have
\begin{displaymath}
\begin{array}{lll}
(G(w)\,,\,w)_X
-\sigma\big(\|u\|^{2}_{V}+\|v\|^{2}_{H}\big)-\frac{\lambda}{2}\|v\|^{2}_{H}
&\geq & (\delta-\sigma)\|u\|_{V}^{2}+(\frac{\lambda}{2}
-\delta-\sigma)\|v\|^{2}_{H}\\\\ & &-\delta\lambda\|u\|_V\|v\|_H\\ \\
&\geq & 2\sqrt{(\delta-\sigma)(\frac{\lambda}{2}
-\delta-\sigma)}\|u\|_V\|v\|_H\\\\ & &-\delta\lambda\|u\|_V\|v\|_H
\end{array}
\end{displaymath}
we can check that $4(\delta-\sigma)(\frac{\lambda}{2}
-\delta-\sigma)=\lambda^2\delta^2$ so that
$$(G(w)\,,\,w)_X
-\sigma\|w\|^{2}_{X}-\frac{\lambda}{2}\|v\|^{2}_{H}\geq 0.$$ The
proof is complete.

The existence of solution for (\ref{wav2.2}) follows in the  same
approach as for equation (\ref{wav3}). Similarly, we can prove an
analogous result as in lemma~\ref{dissipativity} and we have shown
that there also exists a bounded absorbing set $O$ in $X$.

Now let's establish the tail ends estimates  for equation
(\ref{wav2.2}).
\begin{lem}\label{tail2.}
If (\ref{nonlincond1}), (\ref{nonlincond2}) hold, $g\in H$ and
$w_0=(u_0,v_0)\in O$, then for every $\varepsilon >0$, there
exists $T(\varepsilon)$ and $K(\varepsilon)$ such that the
solution $w(t)=(u(t),v(t))$ of problem (\ref{wav2.2}) satisfies
\begin{equation}\label{tail2.0}
\int_{\Omega\cap\{|x|\geq k\}}\Big[|u(t)|^{2} +|\nabla
u(t)|^2+|v(t)|^2\Big]\,dx\leq \varepsilon, \qquad t\geq
T(\varepsilon),\quad k\geq K(\varepsilon).
\end{equation}
\end{lem}
{\bf Proof:} \\The proof works basically like that for equation
(\ref{wav3}). Any solution $w(t)=\left(\begin{array}{ll} u(t)\\
v(t)\end{array}\right)$ satisfies:
\begin{equation}\label{tail2.1}
v_t-\Delta u+ (\lambda -\delta)v+(\delta^2-\lambda\delta)u=-f(u)+g
\end{equation}
and
\begin{equation}\label{tail2.2}
u_t+\delta u=v.
\end{equation}
We choose the same cut-off function $\theta$.\\
Now take inner product in $H$ of $\theta(\frac{|x|^2}{k^2})v(x)$
with (\ref{tail2.1}) to get
\begin{eqnarray}\label{tail2.3}
\begin{array}{lll}\displaystyle
\int_{\Omega}\theta(\frac{|x|^2}{k^2})vv_t\,dx
-\int_{\Omega}\Delta u\theta(\frac{|x|^2}{k^2})v\,dx +
(\lambda-\delta)\int_{\Omega}\theta(\frac{|x|^2}{k^2})|v|^2\,dx
\\
\\\displaystyle +(\delta^2-\lambda\delta)\int_{\Omega}\theta(\frac{|x|^2}{k^2})uv\,dx=
-\int_{\Omega}f(u)\theta(\frac{|x|^2}{k^2})v\,dx+\int_{\Omega}\theta(\frac{|x|^2}{k^2})gv\,dx.
\end{array}
\end{eqnarray}
But
\begin{displaymath}
\begin{array}{lll}\displaystyle
-\int_{\Omega}\Delta u\theta(\frac{|x|^2}{k^2})v\,dx & = &
\displaystyle \int_{\Omega}\theta(\frac{|x|^2}{k^2})\nabla
u\cdot\nabla
v\,\,+\frac{2}{k^2}\int_{\Omega}\theta^\prime(\frac{|x|^2}{k^2})vx\cdot\nabla
u \\ \\
& = & \displaystyle
\int_{\Omega}\theta(\frac{|x|^2}{k^2})\left[\delta|\nabla
u|^2+\nabla u\cdot\nabla u_t
\right]+\frac{2}{k^2}\int_{\Omega}\theta^\prime(\frac{|x|^2}{k^2})vx\cdot\nabla
u \\ \\
&=& \displaystyle
\frac{1}{2}\frac{d}{dt}\int_{\Omega}\theta(\frac{|x|^2}{k^2})|\nabla
u|^2+\delta\int_{\Omega}\theta(\frac{|x|^2}{k^2})|\nabla u|^2\\
\\& &\displaystyle
+\frac{2}{k^2}\int_{\Omega}\theta^\prime(\frac{|x|^2}{k^2})vx\cdot\nabla
u,
\end{array}
\end{displaymath}
and
\begin{displaymath}
\begin{array}{lll}\displaystyle
(\delta^2-\lambda\delta)\int_{\Omega}\theta(\frac{|x|^2}{k^2})uv\,dx
=(\delta^2-\lambda\delta+1)\int_{\Omega}\theta(\frac{|x|^2}{k^2})(\delta
|u|^2+uu_t) \\ \\ = \displaystyle
\frac{1}{2}(\delta^2-\lambda\delta+1)\frac{d}{dt}\int_{\Omega}\theta(\frac{|x|^2}{k^2})
|u|^2+\delta(\delta^2-\lambda\delta+1)\int_{\Omega}\theta(\frac{|x|^2}{k^2})
|u|^2.
\end{array}
\end{displaymath}\\
Then (\ref{tail2.1}) becomes\\
\begin{eqnarray}\label{tail2.4}
\begin{array}{lll}\displaystyle
\frac{1}{2}\frac{d}{dt}\int_{\Omega}\theta(\frac{|x|^2}{k^2})\left[(\delta^2-\lambda\delta)|u|^2+|\nabla
u|^2+|v|^2\right]\\
\\ \displaystyle+\delta\int_{\Omega}\theta(\frac{|x|^2}{k^2})\left[(\delta^2-\lambda\delta)|u|^2+|\nabla
u|^2+|v|^2\right]+
(\lambda-2\delta)\int_{\Omega}\theta(\frac{|x|^2}{k^2})|v|^2\\
\\\displaystyle=-\int_{\Omega}\theta(\frac{|x|^2}{k^2})f(u)(\delta u+u_t)
+\int_{\Omega}\theta(\frac{|x|^2}{k^2})gv\,dx
-\frac{2}{k^2}\int_{\Omega}\theta^\prime(\frac{|x|^2}{k^2})vx\cdot\nabla
u.
\end{array}
\end{eqnarray}\\
But $\delta^2-\lambda\delta$ could be negative for certain  values
of $\lambda$. Since $\delta^2-\lambda\delta+1>0$, let's introduce
another equation  to get a more desirable identity.

Taking inner product of $\theta(\frac{|x|^2}{k^2})u(x)$ with
(\ref{tail2.2}), we get\\
$$\frac{1}{2}\frac{d}{dt}\int_{\Omega}\theta(\frac{|x|^2}{k^2})|u|^2\,dx+\int_{\Omega}\theta(\frac{|x|^2}{k^2})|u|^2\,dx=
\int_{\Omega}\theta(\frac{|x|^2}{k^2})uv\,dx.$$\\
And adding the above and (\ref{tail2.4}) yields\\
\begin{eqnarray}\label{tail2.5}
\begin{array}{lll}\displaystyle
\frac{1}{2}\frac{d}{dt}\int_{\Omega}\theta(\frac{|x|^2}{k^2})\left[(\delta^2-\lambda\delta+1)|u|^2+|\nabla
u|^2+|v|^2\right]\\
\\ \displaystyle+\delta\int_{\Omega}\theta(\frac{|x|^2}{k^2})\left[(\delta^2-\lambda\delta+1)|u|^2+|\nabla
u|^2+|v|^2\right]+
(\lambda-3\delta)\int_{\Omega}\theta(\frac{|x|^2}{k^2})|v|^2\\
\\\displaystyle=-\int_{\Omega}\theta(\frac{|x|^2}{k^2})f(u)(\delta u+u_t)
+\int_{\Omega}\theta(\frac{|x|^2}{k^2})gv\,dx
-\frac{2}{k^2}\int_{\Omega}\theta^\prime(\frac{|x|^2}{k^2})vx\cdot\nabla
u.
\end{array}
\end{eqnarray}\\
Then the conclusion follows the same way as in the proof of
lemma~\ref{tail}.\\

Similarly, we have the following energy equation for the solution
of (\ref{wav2.2}),
\begin{equation}\label{energy.}
\frac{d}{dt}E(w(t))+2\delta E(w(t))=G(w(t))\qquad \forall t>0,
\end{equation}
where
\begin{equation}\label{energy.1}
E(w)=(\delta^2-\lambda\delta+1)\|u\|^{2}_{H}+\|\nabla
u\|^{2}_{H^N}+\|v\|^{2}_{H},
\end{equation}
and
\begin{equation}\label{energy.2}
G(w)=-2(\lambda-3\delta)\|v\|_{H}^{2}+2\int_{\Omega}gv\,dx-2\int_{\Omega}f(u)v\,dx.
\end{equation}

The rest of the proof of existence of a global attractor is again
similar to the case with mass term. We get the main result in this
section
\begin{theo}\label{attractor2}Let $\Omega$ be a domain of $R^N$
bounded in only one direction. Assume that $f$ satisfies
(\ref{nonlincond1}), (\ref{nonlincond2}) and $g\in L^2(\Omega)$.
Then, problem (\ref{wav2.2}) possesses a global attractor in
$X=H^{1}_{0}(\Omega)\times L^2(\Omega)$ which is a compact
invariant subset that attracts every bounded set of $X$ with
respect to the norm topology.
\end{theo}



\bibliographystyle{plain}


{

}\pagebreak  \pagestyle{empty}


\end{document}